\documentclass[a4paper, 12pt]{article}%
\usepackage{amsmath}
\usepackage{graphicx}%
\usepackage{amsfonts}%
\usepackage{amssymb}

\begin{document}

\begin{center}
{\Huge The Szemer\'{e}di property for compact C*-dynamical systems}

\bigskip

Conrad Beyers, Rocco Duvenhage\footnote{{\normalsize Corresponding author.
\textit{E-mail address:} rocco@postino.up.ac.za (R. Duvenhage).}} and Anton Str\"{o}h

\bigskip

\textit{Department of Mathematics and Applied Mathematics}

\textit{University of Pretoria, 0002 Pretoria, South Africa}

\bigskip

2006-5-29
\end{center}

\bigskip

\noindent\textbf{Abstract}

Furstenberg's multiple recurrence result for measure theoretic dynamical
systems is proved for compact $C^{\ast}$-dynamical systems for which the
evolution is given by a semigroup with the right cancellation property, a
right invariant measure and containing a F\o lner net.

\bigskip


\section{Introduction}

In 1977 Furstenberg [2] proved a measure theoretic multiple recurrence theorem
which gave an alternative proof of Szemer\'{e}di's Theorem in combinatorial
number theory. Furstenberg's result states that for a measure preserving
dynamical system $(X,\Sigma,\nu,T)$
\begin{equation}
\liminf_{N\rightarrow\infty}\frac{1}{N}\sum_{n=1}^{N}\nu\left(  V\cap
T^{-n}V\cap T^{-2n}V\cap...\cap T^{-kn}V\right)  >0 \tag{1.1}%
\end{equation}
for any $V\in\Sigma$ with $\nu(V)>0$, where $\nu$ is a probability measure on
the $\sigma$-algebra $\Sigma$ in the set $X$, and $T:X\rightarrow X$ is an
invertible transformation with $\nu\left(  T^{-1}V\right)  =\nu(V)$ for all
$V\in\Sigma$. We view $T^{n}$, where $n$ ranges over $\mathbb{N=}\left\{
1,2,3,...\right\}  $, as the evolution of the system over the semigroup
$\mathbb{N}$. See [3] for a very clear presentation of this result. We will
refer to (1.1) as the Szemer\'{e}di property for measure preserving dynamical
systems. Very roughly put, the proof of Furstenberg's Theorem proceeds by
proving the Szemer\'{e}di property for certain special cases and then properly
combining these cases to prove it in general. One component of the proof is to
consider compact systems, also known as almost periodic systems, namely
systems $(X,\Sigma,\nu,T)$ for which the orbit $\left(  f\circ T^{n}\right)
_{n\in\mathbb{N}}$ for every $f\in L^{2}(\nu)$ is relatively compact
(equivalently, totally bounded) in $L^{2}(\nu)$.

The notion of compactness extends easily to $\ast$-dynamical systems
consisting of a possibly noncommutative $\ast$-algebra $\mathfrak{A}$, a
positive linear functional $\omega$, and an evolution of $\mathfrak{A}$ over a
general semigroup $K$. First steps towards a generalization of Furstenberg's
Theorem to $C^{\ast}$-dynamical systems (i.e. where $\mathfrak{A}$ is a
$C^{\ast}$-algebra) were taken in [5] which included a discussion of compact
systems in the case of an evolution over $\mathbb{N}$. In this paper we prove
the Szemer\'{e}di property for compact $C^{\ast}$-dynamical systems in which
$\omega$ is tracial, i.e. $\omega(ab)=\omega(ba)$ for all $a,b\in\mathfrak{A}%
$, and $K$ is a semigroup with a right invariant measure containing a F\o lner
net (more precisely, $K$ will be a ``F\o lner semigroup'', which we define in
Section 3).

In Section 2 we define compact $\ast$-dynamical systems and obtain a
preliminary recurrence result in terms of seminormed spaces which is used in
Section 4 to prove the Szemer\'{e}di property, Theorem 4.3, for compact
$C^{\ast}$-dynamical systems. F\o lner nets are defined in Section 3, where we
also derive some of their properties which we need in Section 4. Our proof of
the Szemer\'{e}di property in Section 4 follows the basic structure of the one
given in [3], but we have to take into account certain subtleties and
technical difficulties arising from working with a noncommutative $C^{\ast}%
$-algebra rather than with the abelian algebra $L^{\infty}(\nu)$ used in [3],
and with more general groups and semigroups than $\mathbb{Z}$ and $\mathbb{N}%
$. Since we work via abstract seminormed spaces, the structure of the proof
becomes clearer. Essentially the whole paper consists of proving Theorem 4.3
in a series of propositions, lemmas and corollaries (except for Lemma 2.6 and
Proposition 2.7 which are used only in the discussion of an example at the end
of Section 2). One of these, Corollary 4.2, in itself is an interesting
recurrence statement for compact $C^{\ast}$-dynamical systems.

\section{Compact $\ast$-dynamical systems}

In this section we prove a simple recurrence result in seminormed spaces that
has immediate consequences for $\ast$-dynamical systems, which we will define
in a moment.

A linear functional $\omega$ on a $\ast$-algebra $\mathfrak{A}$ is called
positive if $\omega(A^{\ast}A)\geq0$ for all $A\in\mathfrak{A}$. This allows
us to define a seminorm $\Vert\cdot\Vert_{\omega}$ on $\mathfrak{A}$ by
\[
\Vert A\Vert_{\omega}:=\sqrt{\omega(A^{\ast}A)}%
\]
for all $A\in\mathfrak{A}$, as is easily verified using the Cauchy-Schwarz
inequality for positive linear functionals.

\bigskip

\noindent\textbf{Definition 2.1. }Let $\omega$ be a positive linear functional
on a $\ast$-algebra $\mathfrak{A}$, $K$ a semigroup, and $\tau_{g}%
:\mathfrak{A}\rightarrow\mathfrak{A}$ a linear map for each $g\in K$ such
that
\[
\tau_{g}\circ\tau_{h}=\tau_{gh}%
\]
and
\[
\Vert\tau_{g}(A)\Vert_{\omega}=\Vert A\Vert_{\omega}%
\]
for all $g,h\in K$ and $A\in\mathfrak{A}$. Then we call $\left(
\mathfrak{A},\omega,\tau,K\right)  $ a $\ast$\textit{-dynamical system}. If
furthermore $\mathfrak{A}$ is a $C^{\ast}$-algebra and $\Vert\tau_{g}%
(A)\Vert\leq\Vert A\Vert$ in $\mathfrak{A}$'s norm for all $A\in\mathfrak{A}$
and $g\in K$, then we refer to $\left(  \mathfrak{A},\omega,\tau,K\right)  $
as a $C^{\ast}$\textit{-dynamical system}.

\bigskip

Before we proceed, we review definitions and facts that we will need: A set
$V$ in a pseudo metric space $(X,d)$ is said to be $\varepsilon$%
\emph{-separated}, where $\varepsilon>0$, if $d(x,y)\geq\varepsilon$ for all
$x,y\in V$ with $x\neq y$. A set $B\subset X$ is said to be \emph{totally
bounded in} $\left(  X,d\right)  $ if for every $\varepsilon>0$ there exists a
finite set $M_{\varepsilon}\subset X$ such that for every $x\in B$ there is a
$y\in M_{\varepsilon}$ with $d(x,y)<\varepsilon$. It is then not difficult to
show that for any $\varepsilon>0$ there exists a maximal set (in the sense of
cardinality, or number of elements) $V\subset B$ that is $\varepsilon
$-separated, and furthermore, if $B\neq\varnothing$, then $V$ is finite with
$|V|>0$.

\bigskip

\noindent\textbf{Definition 2.2.} A $\ast$-dynamical system $\left(
\mathfrak{A},\omega,\tau,K\right)  $ is called \textit{compact} if the
\textit{orbit}
\[
B_{A}:=\{\tau_{g}(A):g\in K\}
\]
is totally bounded in $(\mathfrak{A},\Vert\cdot\Vert_{\omega})$ for each
$A\in\mathfrak{A}$.

\bigskip

\noindent\textbf{Definition 2.3.} Let $K$ be a semigroup. We call a set
$E\subset K$ \emph{relatively dense in }$K$ if there exist an $r\in\mathbb{N}$
and $g_{1},\ldots,g_{r}\in K$ such that
\[
E\cap\{gg_{1},\ldots,gg_{r}\}\neq\varnothing
\]
for all $g\in K$.

\bigskip

Strictly speaking one could call this \textit{left} relatively denseness, with
the right hand case being defined similarly in terms of $g_{j}g$, but we will
only work with Definition 2.3 in this paper. The usual definition of relative
denseness of a subset $E$ in $\mathbb{N}$ is in terms of ``bounded gaps'' (see
[6] for example), and it is easy to check that in this special case the two
definitions are equivalent.

\bigskip

\noindent\textbf{Proposition 2.4. }\textit{Let }$K$\textit{ be a semigroup,
}$\left(  X,\Vert\cdot\Vert\right)  $\textit{ a seminormed space, and }%
$U_{g}:X\rightarrow X$\textit{ a linear map for each }$g\in K$\textit{ such
that }$U_{g}U_{h}=U_{gh}$\textit{ and }$\Vert U_{g}x\Vert\geq\Vert x\Vert
$\textit{ for all }$g,h\in K$\textit{ and }$x\in X$\textit{. Suppose that}
$B_{x_{0}}:=\{U_{g}x_{0}:g\in K\}$\textit{ is totally bounded in }$\left(
X,\Vert\cdot\Vert\right)  $\textit{ for some }$x_{0}\in X$\textit{. }%
\noindent\textit{Then for each }$\varepsilon>0$\textit{, the set}
\[
E:=\{g\in K:\Vert U_{g}x_{0}-x_{0}\Vert<\varepsilon\}
\]
\textit{is relatively dense in }$K$.

\bigskip

\noindent\textbf{Proof. }Since $B_{x_{0}}$ is totally bounded in $\left(
X,\Vert\cdot\Vert\right)  $, there is a maximal $V=\{U_{g_{1}}x_{0}%
,...,U_{g_{r}}x_{0}\}$, with $U_{g_{j}}x_{0}\neq U_{g_{l}}x_{0}$ whenever
$j\neq l$, which is $\varepsilon$-separated. \noindent But $\Vert
U_{g^{\prime}gg_{j}}x_{0}-U_{g^{\prime}gg_{l}}x_{0}\Vert\geq\Vert U_{g_{j}%
}x_{0}-U_{g_{l}}x_{0}\Vert$ for any $g,g^{\prime}\in K$, hence $V_{g^{\prime
}g}:=\{U_{g^{\prime}gg_{1}}x_{0},...,U_{g^{\prime}gg_{r}}x_{0}\}$ is
$\varepsilon$-separated, with $r$ elements. Since $V_{g^{\prime}g}\subset
B_{x_{0}}$, it is also maximally $\varepsilon$-separated in $B_{x_{0}}$. But
$U_{g^{\prime}}x_{0}\in B_{x_{0}}$, therefore $\Vert U_{gg_{j}}x_{0}%
-x_{0}\Vert\leq\Vert U_{g^{\prime}gg_{j}}x_{0}-U_{g^{\prime}}x_{0}%
\Vert<\varepsilon$ for some $j\in\{1,\ldots,r\}$. Hence, for each $g\in K$
there exists an $h\in\{gg_{1},\ldots,gg_{r}\}$ such that $\Vert U_{h}%
x_{0}-x_{0}\Vert<\varepsilon$, i.e.
\[
E\cap\{gg_{1},\ldots,gg_{r}\}\neq\varnothing
\]
for all $g\in K$, and so $E$ is relatively dense in $K$. $\square$

\bigskip

\noindent\textbf{Corollary 2.5.} \textit{Let }$\left(  \mathfrak{A}%
,\omega,\tau,K\right)  $\textit{ be a compact }$\ast$\textit{-dynamical system
and let }$m_{0},...,m_{k}\in\mathbb{N}\cup\{0\}$\textit{. For any
}$\varepsilon>0$\textit{ and }$A\in\mathfrak{A}$\textit{, the set}
\[
E:=\{g\in K:\Vert\tau_{g^{m_{j}}}(A)-A\Vert_{\omega}<\varepsilon\text{ for
}j=0,...,k\}
\]
\textit{is then relatively dense in }$K$\textit{, where we write }$\tau
_{g^{0}}(A)\equiv A$.

\bigskip

\noindent\textbf{Proof.} Without loss we can assume that none of the $m_{j}$
's are zero. Then the result follows from Proposition 2.4 with $\varepsilon$
replaced by $\varepsilon/\max\{m_{0},\ldots,m_{k}\}$, since for every
$j=0,...,k$ we have
\begin{align*}
&  \Vert\tau_{g^{m_{j}}}(A)-A\Vert_{\omega}\\
&  \leq\Vert\tau_{g^{m_{j}}}(A)-\tau_{g^{m_{j}-1}}(A)\Vert_{\omega}+\Vert
\tau_{g^{m_{j}-1}}(A)-\tau_{g^{m_{j}-2}}(A)\Vert_{\omega}+\ldots+\Vert\tau
_{g}(A)-A\Vert_{\omega}\\
&  =\Vert\tau_{g^{m_{j}-1}}[\tau_{g}(A)-A]\Vert_{\omega}+\Vert\tau
_{g^{m_{j}-2}}[\tau_{g}(A)-A]\Vert_{\omega}+\ldots+\Vert\tau_{g}%
(A)-A\Vert_{\omega}\\
&  =m_{j}\Vert\tau_{g}(A)-A\Vert_{\omega}\\
&  <\varepsilon
\end{align*}
for all $g\in K$ for which $\Vert\tau_{g}(A)-A\Vert_{\omega}<\varepsilon
/\max\{m_{0},\ldots,m_{k}\}$. $\square$

\bigskip

We now briefly indicate why the $L^{2}$ definition of compactness given in
Section 1 is a special case of Definition 2.2. Given a $\ast$-dynamical system
$\left(  \mathfrak{A},\omega,\tau,K\right)  $, the GNS construction provides
us with a representation of $\left(  \mathfrak{A},\omega\right)  $, namely an
inner product space $\mathfrak{G}$, a linear surjection $\iota:\mathfrak{A}%
\rightarrow\mathfrak{G}$, and a linear mapping $\pi:\mathfrak{A}\rightarrow
L(\mathfrak{G})$, with $L(\mathfrak{G})$ the space of all linear maps
$\mathfrak{G}\rightarrow\mathfrak{G}$ (not necessarily bounded), such that
$\left\langle \iota(A),\iota(B)\right\rangle =\omega(A^{\ast}B)$, $\pi
(A)\iota(B)=\iota(AB)$ and $\pi(AB)=\pi(A)\pi(B)$ for all $A,B\in\mathfrak{A}%
$. Then
\[
U_{g}:\mathfrak{G}\rightarrow\mathfrak{G}:\iota(A)\mapsto\iota(\tau_{g}(A))
\]
is a well-defined linear operator with $\left|  \left|  U_{g}x\right|
\right|  =\left|  \left|  x\right|  \right|  $ for all $x\in\mathfrak{G}$ and
$g\in K$. It is then straightforward to show that $\left(  \mathfrak{A}%
,\omega,\tau,K\right)  $ is compact if and only if all the orbits
\begin{equation}
B_{x}:=\left\{  U_{g}x:g\in K\right\}  \tag{2.1}%
\end{equation}
with $x\in\mathfrak{G}$, are totally bounded in $\mathfrak{G}$. However,
$U_{g}$ has a unique continuous extension to the completion $\mathfrak{H}$ of
$\mathfrak{G}$, and one can show that all the orbits $B_{x}$, $x\in
\mathfrak{H}$, again defined as in (2.1), are totally bounded if and only if
they are totally bounded for all $x\in\mathfrak{G}$. Hence $\left(
\mathfrak{A},\omega,\tau,K\right)  $ is compact if and only if all the orbits
$B_{x}$, $x\in\mathfrak{H}$, in the Hilbert space $\mathfrak{H}$ are totally
bounded. The measure theoretic definition in Section 1 is a special case of
this simply because $L^{2}(\nu)$ is a Hilbert space obtained exactly as
$\mathfrak{H}$ above through the GNS-construction applied to the state
$\omega=\int(\cdot)d\nu$ on the $C^{\ast}$-algebra $B_{\infty}(\Sigma)$ of all
bounded complex-valued $\Sigma$-measurable functions on $X$, or to $\omega$ on
$L^{\infty}(\nu)$.

To conclude this section we present an example of a compact $C^{\ast}%
$-dynamical system in which the $C^{\ast}$-algebra is noncommutative. To do
this we need a few simple tools, which we now discuss.

First note that if a set in a $C^{\ast}$-algebra $\mathfrak{A}$ is totally
bounded in $\mathfrak{A}$ (i.e. in terms of $\mathfrak{A}$'s norm), then it is
also totally bounded in $\left(  \mathfrak{A},\left|  \left|  \cdot\right|
\right|  _{\omega}\right)  $ for any positive linear functional $\omega$ on
$\mathfrak{A}$, since $\left|  \left|  \cdot\right|  \right|  _{\omega}%
\leq\left|  \left|  \omega\right|  \right|  ^{1/2}\left|  \left|
\cdot\right|  \right|  $ as is easily verified (keep in mind that $\omega$ is
bounded, since it is positive and $\mathfrak{A}$ is a $C^{\ast}$-algebra).
Hence, if we can prove that the orbits of a given $C^{\ast}$-dynamical system
$\left(  \mathfrak{A},\omega,\tau,K\right)  $ are totally bounded in
$\mathfrak{A}$, then it follows that the system is compact. Of course, this is
then a stronger form of compactness, but Example 2.8 happens to possess this
stronger property, and it turns out to be easier to prove than to prove
compactness directly in terms of $\left|  \left|  \cdot\right|  \right|
_{\omega}$, since $\mathfrak{A}$'s norm is submultiplicative, which makes it
easier to work with than $\left|  \left|  \cdot\right|  \right|  _{\omega}$.

In the remainder of this section we work with a $C^{\ast}$-algebra
$\mathfrak{A}$, an arbitrary \textit{set} $K$, and a $\ast$-homomorphism
$\tau_{g}:\mathfrak{A}\rightarrow\mathfrak{A}$ for each $g\in K$. So for the
moment we are working in a more general setting than Definition 2.1, however
in Example 2.8 below we will be more specific. When we say that an ``orbit''
$(a_{g})\equiv\left(  a_{n}\right)  _{g\in K}$ is \textit{totally bounded} in
a space, we mean that the set $\{a_{g}:g\in K\}$ is totally bounded in that
space. For any subset $\mathfrak{V}\subset\mathfrak{A}$ we will denote the set
of all polynomials over $\mathbb{C}$ generated by the elements of
$\mathfrak{V}$ and their adjoints, by $p(\mathfrak{V})$, i.e. $p(\mathfrak{A}%
)$ consists of all finite linear combinations of all finite products of
elements of $\mathfrak{V}\cup\mathfrak{V}^{\ast}$ with $\mathfrak{V}^{\ast
}:=\left\{  a^{\ast}:a\in\mathfrak{V}\right\}  $. We will use the notation
$XY:=\left\{  xy:x\in X,y\in Y\right\}  $ whenever $X$ and $Y$ are sets for
which this multiplication of their elements is defined.

\bigskip

\noindent\textbf{Lemma 2.6.}\textit{ If }$\left(  \tau_{g}(A)\right)
$\textit{ is totally bounded in $\mathfrak{A}$ for every }$A$\textit{ in some
subset $\mathfrak{V}$ of $\mathfrak{A}$, then }$\left(  \tau_{g}(A)\right)
$\textit{ is totally bounded in $\mathfrak{A}$ for every }$A\in p(\mathfrak{V}%
)$\textit{.}

\bigskip

\noindent\textbf{Proof.} Consider any $A,B\in\mathfrak{A}$ for which $\left(
\tau_{g}(A)\right)  $ and $\left(  \tau_{g}(B)\right)  $ are totally bounded
in $\mathfrak{A}$, and any $\varepsilon>0$. By the hypothesis there are finite
sets $M,N\subset\mathfrak{A}$ such that for each $g\in K$ there is an
$a_{g}\in M$ and a $b_{g}\in N$ such that $\Vert\tau_{g}(A)-a_{g}%
\Vert<\varepsilon$ and $\Vert\tau_{g}(B)-b_{g}\Vert<\varepsilon$. Clearly
\begin{align*}
\Vert\tau_{g}(A)\tau_{g}(B)-a_{g}b_{g}\Vert &  \leq\Vert\tau_{g}(A)\Vert
\Vert\tau_{g}(B)-b_{g}\Vert+\Vert\tau_{g}(A)-a_{g}\Vert\Vert b_{g}\Vert\\
&  \leq\varepsilon\left(  \Vert\tau_{g}(A)\Vert+\Vert b_{g}\Vert\right)
\end{align*}
but note that $\Vert\tau_{g}(A)\Vert\leq\Vert A\Vert$, since $\tau_{g}$ is a
$\ast$-homomorphism and $\mathfrak{A}$ is a $C^{\ast}$-algebra, while $\left|
\left|  b_{g}\right|  \right|  <\Vert\tau_{g}(B)\Vert+\varepsilon\leq\left|
\left|  B\right|  \right|  +\varepsilon$. Since $MN$ is a finite subset of
$\mathfrak{A}$, and $a_{g}b_{g}\in MN$, it follows that $(\tau_{g}(AB))$ is
totally bounded in $\mathfrak{A}$. Similarly $(\tau_{g}(A^{\ast}))$ and
$(\tau_{g}(\alpha A+\beta B))$ are totally bounded in $\mathfrak{A}$ for any
$\alpha,\beta\in\mathbb{C}$, and this is enough to prove the lemma. $\square$

\bigskip

\noindent\textbf{Proposition 2.7.}\textit{ Now assume that }$\mathfrak{A}%
$\textit{ is generated by a subset }$\mathfrak{V}\subset\mathfrak{A}$\textit{
for which }$\tau_{g}(\mathfrak{V})\subset p(\mathfrak{V})$\textit{ for every
}$g\in K$\textit{. Also assume that }$(\tau_{g}(A))$\textit{ is totally
bounded in $\mathfrak{A}$ for every }$A\in\mathfrak{V}$\textit{. Then }%
$(\tau_{g}(A))$\textit{ is totally bounded in $\mathfrak{A}$ for every }%
$A\in\mathfrak{A}$\textit{.}

\bigskip

\noindent\textbf{Proof.} Firstly it is easily shown that if $Y$ is a dense
subspace of a normed space $X$, $U_{g}:Y\rightarrow Y$ is linear with $\Vert
U_{g}\Vert\leq1$ for all $g\in K$, and $\left(  U_{g}y\right)  $ is totally
bounded in $X$ for every $y\in Y$ (or in $Y$ for every $y\in Y$), then for the
unique bounded linear extension $U_{g}:X\rightarrow X$ the ``orbit'' $\left(
U_{g}x\right)  $ is totally bounded in $X$ for every $x\in X$. (We also used
this fact when we discussed the GNS-construction above.)

Now simply set $X=\mathfrak{A}$, $Y=p(\mathfrak{V})$ and $U_{g}=\tau_{g}$,
then by our assumptions and Lemma 2.6 all the requirements in the remark above
are met. $\square$

\bigskip

\noindent\textbf{Example 2.8.} We consider a so-called rotation $C^{\ast}%
$-algebra, and use Proposition 2.7 to show that we obtain a compact $C^{\ast}%
$-dynamical system. As described in [1], let $\mathfrak{H}:=L^{2}%
(\mathbb{R}/\mathbb{Z})$ and define two unitary operators $U$ and $V$ on
$\mathfrak{H}$ by
\[
\left(  Uf\right)  (t)=f(t+\theta)
\]
and
\[
\left(  Vf\right)  (t)=e^{2\pi it}f(t)
\]
for $f\in\mathfrak{H}$, where $\theta\in\mathbb{R}$ (though the interesting
case is $\theta\in\mathbb{Q}$). These operators satisfy
\begin{equation}
UV=e^{2\pi i\theta}VU\text{.} \tag{2.2}%
\end{equation}
Let $\mathfrak{A}$ be the $C^{\ast}$-algebra generated by $U$ and $V$. Note
that $\mathfrak{A}$ is noncommutative because of (2.2). Then, as shown in [1],
there is a unique trace $\omega$ on $\mathfrak{A}$, i.e. a state with
$\omega(AB)=\omega(BA)$ (we will return to traces in Section 4). Define
$\tau:\mathfrak{A}\rightarrow\mathfrak{A}$ by $\tau(A)=U^{\ast}AU$ for all
$A\in\mathfrak{A}$, then $\tau$ is a $\ast$-isomorphism and therefore $\left|
\left|  \tau(A)\right|  \right|  =\left|  \left|  A\right|  \right|  $, since
$\mathfrak{A}$ is a $C^{\ast}$-algebra. Also, since $\omega$ is a trace and
$U$ is unitary, $\Vert\tau(A)\Vert_{\omega}=\Vert A\Vert_{\omega}$ for all
$A\in\mathfrak{A}$. Hence $(\mathfrak{A},\omega,\tau,\mathbb{N})$ is a
$C^{\ast}$-dynamical system, where by slight abuse of notation $\tau$ here
denotes the function $n\mapsto\tau^{n}$ as well, to fit it into Definition
2.1's notation.

We now show that $(\mathfrak{A},\omega,\tau,\mathbb{N})$ is compact: It is
trivial that $(\tau^{n}(U))=(U)$ is totally bounded in $\mathfrak{A}$.
Furthermore, $\tau^{n}(V)=(U^{\ast})^{n}VU^{n}=e^{-2\pi in\theta}V$ by (2.2).
Since the unit circle is compact, it follows that $(\tau^{n}(V))$ is totally
bounded in $\mathfrak{A}$. From Proposition 2.7 with $\mathfrak{V}=\left\{
U,V\right\}  $ we conclude that $(\tau^{n}(A))$ is totally bounded in
$\mathfrak{A}$ for all $A\in\mathfrak{A}$. In particular the $C^{\ast}%
$-dynamical system $(\mathfrak{A},\omega,\tau,\mathbb{N})$ is compact.
Similarly $(\mathfrak{A},\omega,\tau,\mathbb{Z})$ is compact.

\section{F\o lner nets}

In this section we define F\o lner nets in an abstract setting (see [4] for
some discussion in the more specific case of topological groups) and present a
number of facts regarding these nets leading to Proposition 3.6 which we will
need to prove our main result, Theorem 4.3, in the next section.

To simplify statements of definitions and results in the sequel we introduce
the following terminology: A triple $\left(  K,\Sigma,\mu\right)  $ with $K$ a
semigroup, $\Sigma$ a $\sigma$-algebra in $K$, and $\mu$ a positive measure on
$\Sigma$, will be called a \textit{measure semigroup}. When $\mu$ is right
invariant, i.e. $Vg\in\Sigma$ and $\mu(Vg)=\mu(V)$ for $V\in\Sigma$ and $g\in
K$, we say that $\left(  K,\Sigma,\mu\right)  $ is \textit{right invariant}.
When we say that a net $\left(  \Lambda_{\alpha}\right)  $ has some property
for $\alpha$ ``large enough'', then we mean that there is a $\beta$ in the
directed set such that the property holds for all $\alpha\geq\beta$.

\bigskip

\noindent\textbf{Definition 3.1.} Let $\left(  K,\Sigma,\mu\right)  $ be a
measure semigroup such that $g\Lambda\in\Sigma$ for all $g\in K$ and
$\Lambda\in\Sigma$. A net $\left(  \Lambda_{\alpha}\right)  $ in $\Sigma$ is
called a \textit{F\o lner net} in $\left(  K,\Sigma,\mu\right)  $ if
$0<\mu(\Lambda_{\alpha})<\infty$ for $\alpha$ large enough and
\[
\lim_{\alpha}\frac{\mu\left(  \Lambda_{\alpha}\Delta(g\Lambda_{\alpha
})\right)  }{\mu(\Lambda_{\alpha})}=0
\]
for all $g\in K$. When these conditions are satisfied, in particular such a
net $\left(  \Lambda_{\alpha}\right)  $ exists, and furthermore $K$ has the
right cancellation property (i.e. $g_{1}h=g_{2}h\Rightarrow g_{1}=g_{2}$),
$\mu$ is right invariant, and lastly, if $V\subset K$ with $Vg\in\Sigma$ for
some $g\in K$ implies that $V\in\Sigma$, then we will call $\left(
K,\Sigma,\mu\right)  $ a \textit{F\o lner semigroup}.

\bigskip

Simple examples of (abelian) F\o lner semigroups are $\mathbb{N}$ with the
counting measure, and $[0,\infty)$, the first quadrant in $\mathbb{R}^{2}$,
etc., with Lebesgue measure.

\bigskip

\noindent\textbf{Lemma 3.2.}\textit{ Consider a right invariant measure
semigroup }$\left(  K,\Sigma,\mu\right)  $\textit{ which has the right
cancellation property and contains a F\o lner net }$\left(  \Lambda_{\alpha
}\right)  $\textit{. Take any }$g_{\alpha}\in K$\textit{ for each }$\alpha
$\textit{ in the directed set of the net. Then the net}
\[
\left(  \Lambda_{\alpha}g_{\alpha}\right)
\]
\textit{is also a F\o lner net in }$\left(  K,\Sigma,\mu\right)  $.

\bigskip

\noindent\textbf{Proof.} Since $K$ has the right cancellation property, we
have $(Ag)\Delta(Bg)=(A\Delta B)g$ for all $A,B\subset K$ and $g\in
K$.\ Hence
\begin{align*}
\frac{\mu\left(  (\Lambda_{\alpha}g_{\alpha})\Delta(g(\Lambda_{\alpha
}g_{\alpha}))\right)  }{\mu(\Lambda_{\alpha}g_{\alpha})}  &  =\frac{\mu\left(
(\Lambda_{\alpha}\Delta(g\Lambda_{\alpha}))g_{\alpha}\right)  }{\mu
(\Lambda_{\alpha}g_{\alpha})}\\
&  =\frac{\mu\left(  \Lambda_{\alpha}\Delta(g\Lambda_{\alpha})\right)  }%
{\mu(\Lambda_{\alpha})}\\
&  \longrightarrow0
\end{align*}
with respect to $\alpha$. $\square$

\bigskip

\noindent\noindent\textbf{Definition 3.3. }Let $\left(  K,\Sigma,\mu\right)  $
be a measure semigroup. Let $\left(  \Lambda_{\alpha}\right)  $ be a net in
$\Sigma$ with $0<\mu(\Lambda_{\alpha})<\infty$ for $\alpha$ large enough.
Consider any $V\in\Sigma$ and set
\[
D_{\left(  \Lambda_{\alpha}\right)  }(V):=\lim_{\alpha}\left[  \inf\left\{
\frac{\mu(\Lambda_{\beta}\cap V)}{\mu(\Lambda_{\beta})}:\beta\geq
\alpha\right\}  \right]  \equiv\liminf_{\alpha}\frac{\mu(\Lambda_{\alpha}\cap
V)}{\mu(\Lambda_{\alpha})}.
\]
If $D_{\left(  \Lambda_{\alpha}\right)  }(V)>0$, then we say that\ $V$ has
\emph{positive lower density relative to }$\left(  \Lambda_{\alpha}\right)  $.

\bigskip

It is easily checked that $D_{\left(  \Lambda_{\alpha}\right)  }(V)$ in this
definition always exists.

\bigskip

\noindent\textbf{Lemma 3.4. }\textit{Let }$\left(  K,\Sigma,\mu\right)
$\textit{ be a right invariant measure semigroup. Assume that if }$V\subset
K$\textit{ and }$Vg\in\Sigma$\textit{ for some }$g\in K$\textit{, then }%
$V\in\Sigma$\textit{. Let }$E\in\Sigma$\textit{ be relatively dense in }%
$K$\textit{. Then there exists an }$r\in\mathbb{N}$\textit{ and }$g_{1}%
,\ldots,g_{r}\in K$\textit{ such that the following holds: for each }%
$B\in\Sigma$\textit{ with }$\mu(B)<\infty$\textit{ there exists a }%
$j\in\{1,\ldots,r\}$\textit{ such that}%
\[
\mu((Bg_{j})\cap E)\geq\frac{1}{r}\mu(B)\text{.}%
\]
\textbf{Proof.}\ Let $g_{1},...,g_{r}$ be given by Definition 2.3. Set
$B_{j}:=\{b\in B:bg_{j}\in E\}$ for $j=1,\ldots,r$, so $B_{j}g_{j}%
=(Bg_{j})\cap E\in\Sigma$ and hence $B_{j}\in\Sigma$. Now, for any $b\in B$ we
know from Definition 2.3 that $E\cap\{bg_{1},\ldots,bg_{r}\}\neq\varnothing$.
So $bg_{j}\in E$ for some $j\in\{1,\ldots,r\}$, i.e. $b\in B_{j}$. Hence
$B=\bigcup_{j=1}^{r}B_{j}$ and therefore
\[
\mu(B)=\mu(\bigcup_{j=1}^{r}B_{j})\leq\sum_{j=1}^{r}\mu(B_{j})=\sum_{j=1}%
^{r}\mu(B_{j}g_{j})=\sum_{j=1}^{r}\mu((Bg_{j})\cap E)
\]
from which the conclusion follows. $\square$

\bigskip

\noindent\textbf{Lemma 3.5. }\textit{Let }$\left(  K,\Sigma,\mu\right)
$\textit{ be a F\o lner semigroup. \noindent Let }$E\in\Sigma$\textit{ be
relatively dense in }$K$\textit{. Then }$E$\textit{ has positive lower density
relative to some F\o lner net in }$\left(  K,\Sigma,\mu\right)  $\textit{.}

\bigskip

\noindent\textbf{Proof. }Consider any F\o lner net $\left(  \Lambda_{\alpha
}\right)  $ in $\left(  K,\Sigma,\mu\right)  $, then for $\alpha$ large enough
we have $0<\mu(\Lambda_{\alpha})<\infty$ and we will now work only with such
$\alpha$'s without loss of generality. Let $g_{1},\ldots,g_{r}\in K$ be as in
Definition 2.3. \noindent For each $\alpha$ it follows from Lemma 3.4 that
there exists a $j(\alpha)\in\{1,\ldots,r\}$ such that
\[
\frac{\mu((\Lambda_{\alpha}g_{j(\alpha)})\cap E)}{\mu(\Lambda_{\alpha
}g_{j(\alpha)})}\geq\frac{1}{r}%
\]
where we also made use of $\mu(\Lambda_{\alpha}g_{j(\alpha)})=\mu
(\Lambda_{\alpha})$. But it follows from Lemma 3.2 that $\left(
\Lambda_{\alpha}^{\prime}\right)  $ given by $\Lambda_{\alpha}^{\prime
}:=\Lambda_{\alpha}g_{j(\alpha)}$ is a F\o lner net in $\left(  K,\Sigma
,\mu\right)  $. Furthermore,
\begin{align*}
D_{\left(  \Lambda_{\alpha}^{\prime}\right)  }(E)  &  =\liminf_{\alpha
}\frac{\mu(\Lambda_{\alpha}^{\prime}\cap E)}{\mu(\Lambda_{\alpha}^{\prime})}\\
&  =\lim_{\alpha}\left[  \inf\left\{  \frac{\mu(\Lambda_{\beta}^{\prime}\cap
E)}{\mu(\Lambda_{\beta}^{\prime})}:\beta\geq\alpha\right\}  \right] \\
&  \geq\lim_{\alpha}\frac{1}{r}\\
&  =\frac{1}{r}\text{. }\square
\end{align*}

\noindent\textbf{Proposition 3.6. }\textit{Let }$\left(  K,\Sigma,\mu\right)
$\textit{ be a F\o lner semigroup, and }$f:K\rightarrow\mathbb{R}$\textit{ a
}$\Sigma$\textit{-measurable function with }$f\geq0$\textit{.} \textit{Assume
that }$f(g)\geq a$\textit{ for some }$a>0$\textit{ and all }$g$\textit{ in
some relatively dense }$E\in\Sigma$\textit{ in }$K$\textit{. Then there exists
a F\o lner net }$(\Lambda_{\alpha})$\textit{ in }$\left(  K,\Sigma,\mu\right)
$\textit{ such that}
\[
\liminf_{\alpha}\frac{1}{\mu(\Lambda_{\alpha})}\int_{\Lambda_{\alpha}}%
fd\mu>0\text{.}%
\]

\noindent\textbf{Proof. }By Lemma 3.5 there exists a F\o lner net
$(\Lambda_{\alpha})$ in $\left(  K,\Sigma,\mu\right)  $ such that
\begin{align*}
\liminf_{\alpha}\frac{1}{\mu(\Lambda_{\alpha})}\int_{\Lambda_{\alpha}}fd\mu &
\geq\liminf_{\alpha}\frac{1}{\mu(\Lambda_{\alpha})}\int_{\Lambda_{\alpha}\cap
E}fd\mu\\
&  \geq\liminf_{\alpha}\frac{1}{\mu(\Lambda_{\alpha})}\int_{\Lambda_{\alpha
}\cap E}a~dg\\
&  =aD_{\left(  \Lambda_{\alpha}\right)  }(E)\\
&  >0\text{.}~\square
\end{align*}

\section{The Szemer\'{e}di property}

We now put the work of the previous two sections together to prove the
Szemer\'{e}di property for compact $C^{\ast}$-dynamical systems as defined in
Definition 2.1. The first two results below are independent of F\o lner nets
and the results of Section 3, and hence hold in a more general context than
the third and final result, Theorem 4.3 (the Szemer\'{e}di property), which
rests on all that went before it, except Lemma 2.6 and Proposition 2.7.

We will work with a $C^{\ast}$-algebra $\mathfrak{A}$ and take $\omega$ to be
a positive linear functional on $\mathfrak{A}$. It then follows that $\omega$
is bounded, and without loss we can assume that $\left|  \left|
\omega\right|  \right|  =1$ (the case $\omega=0$ being trivial), i.e. $\omega$
is a \textit{state} on $\mathfrak{A}$. By the Cauchy-Schwarz inequality we
have
\[
\left|  \omega(AB)\right|  \leq\left|  \left|  A^{\ast}\right|  \right|
_{\omega}\left|  \left|  B\right|  \right|  _{\omega}\leq\sqrt{\left|  \left|
AA^{\ast}\right|  \right|  }\left|  \left|  B\right|  \right|  _{\omega
}=\left|  \left|  A\right|  \right|  \left|  \left|  B\right|  \right|
_{\omega}%
\]
A \textit{trace} is defined to be a state $\omega$ on a $C^{\ast}$-algebra
$\mathfrak{A}$ such that $\omega(AB)=\omega(BA)$ for all $A,B\in\mathfrak{A}$.
Note that from the previous inequality we then we have
\[
|\omega(ABC)|=|\omega(CAB)|\leq\Vert A\Vert\Vert B\Vert_{\omega}\Vert C\Vert
\]
for all $A,B,C\in\mathfrak{A}$. This fact is used in the proof of Proposition
4.1, along with the following identity which holds in any algebra
$\mathfrak{A}$ and is easily verified by induction:%
\[
\prod_{j=0}^{k}a_{j}-\prod_{j=0}^{k}b_{j}=\sum_{j=0}^{k}\left(  \prod
_{l=0}^{j-1}a_{l}\right)  \left(  a_{j}-b_{j}\right)  \left(  \prod
_{l=j+1}^{k}b_{l}\right)
\]
for any $a_{j},b_{j}\in\mathfrak{A}$.

\bigskip

\noindent\textbf{Lemma 4.1.} \textit{Let }$\mathfrak{A}$\textit{ be a
}$C^{\ast}$\textit{-algebra and }$\omega$\textit{ a trace on }$\mathfrak{A}%
$\textit{. Suppose that }$b\in A^{+}$\textit{, }$\Vert b\Vert\leq1$\textit{
and }$\omega(b)>0$\textit{. Let }$k\in N\cup\{0\}$\textit{, then}\emph{
}$\omega(b^{k+1})>0$\textit{ so we can choose }$\varepsilon>0$\textit{ such
that }$\varepsilon<\omega(b^{k+1})/(k+1)$\textit{. \noindent Set }%
$a:=\omega(b^{k+1})-(k+1)\varepsilon$. \textit{Consider }$c_{0},\ldots
,c_{k}\in\mathfrak{A}$\textit{ such that }$\Vert c_{j}\Vert\leq1$\textit{ and
}$\Vert c_{j}-b\Vert_{\omega}<\varepsilon$\emph{ }\textit{for}$~j=0,\ldots
,k$\textit{. Then}
\[
\left|  \omega\left(  \prod_{j=0}^{k}c_{j}\right)  \right|  >a>0\text{.}%
\]

\noindent\textbf{Proof.} We have $\omega(b^{k+1})>0$ by using the Gelfand
representation of the abelian $C^{\ast}$-algebra $\mathfrak{B}$ generated by
$b$, restricting $\omega$ to $\mathfrak{B}$ and then using Riesz's theorem to
represent $\omega$ by a positive measure on the locally compact Hausdorff
space appearing in the Gelfand representation. \noindent Furthermore,
\begin{align*}
\left|  \omega\left(  \prod_{j=0}^{k}c_{j}\right)  -\omega(b^{k+1})\right|
&  =\left|  \omega\left(  \prod_{j=0}^{k}c_{j}-\prod_{j=0}^{k}b\right)
\right| \\
&  =\left|  \omega\left(  \sum_{j=0}^{k}\left(  \prod_{l=0}^{j-1}c_{l}\right)
\left(  c_{j}-b\right)  \left(  \prod_{l=j+1}^{k}b\right)  \right)  \right| \\
&  \leq\sum_{j=0}^{k}\left(  \left\|  \prod_{l=0}^{j-1}c_{l}\right\|  \left\|
c_{j}-b\right\|  _{\omega}\left\|  b^{k-j}\right\|  \right) \\
&  \leq\sum_{j=0}^{k}\left\|  c_{j}-b\right\|  _{\omega}\\
&  <(k+1)\varepsilon.
\end{align*}
Hence
\[
\left|  \omega\left(  \prod_{j=0}^{k}c_{j}\right)  \right|  >\omega
(b^{k+1})-(k+1)\varepsilon=a>0\text{. }\square
\]

\bigskip

\noindent\textbf{Corollary 4.2.} \emph{Let }$(\mathfrak{A},\omega,\tau
,K)$\emph{ be a compact $C^{\ast}$-dynamical system with $\omega$ a trace.
Suppose that $A\in\mathfrak{A}^{+},$ and $\omega(A)>0$. Take any $m_{0}%
,\ldots,m_{k}\in\mathbb{N}\cup\{0\}$.} \noindent\emph{Then there exists a
relatively dense set $E$ in }$K$\emph{ and an $a>0$ such that }%
\begin{equation}
\left|  \omega\left(  \prod_{j=0}^{k}\tau_{g^{m_{j}}}(A)\right)  \right|  >a
\tag{4.1}%
\end{equation}
\emph{for all $g\in E$.}

\bigskip

\noindent\textbf{Proof.} Since $\omega(A)>0$, $\left|  \left|  A\right|
\right|  >0$, so we can set $b:=A/\left|  \left|  A\right|  \right|  $. For
$c_{j}:=\tau_{g^{m_{j}}}(b)$ we have $\left|  \left|  c_{j}\right|  \right|
\leq\left|  \left|  b\right|  \right|  =1$, so from Lemma 4.1 it follows that
there is an $a^{\prime}>0$ such that
\[
\left|  \omega\left(  \prod_{j=0}^{k}\tau_{g^{m_{j}}}(b)\right)  \right|
>a^{\prime}%
\]
for every $g\in K$ for which $\Vert\tau_{g^{m_{j}}}(b)-b\Vert_{\omega
}<\varepsilon<\omega(b^{k+1})/(k+1)$ for all $j=0,...,k$. By Corollary 2.5
this set of $g$ 's is relatively dense in $K$. Now simply set $a=a^{\prime
}\left|  \left|  A\right|  \right|  ^{k+1}$. $\square$

\bigskip

Finally we reach our set goal, namely a Szemer\'{e}di property for compact
$C^{\ast}$-dynamical systems, which together with Corollary 4.2 form the main
results of this paper:

\bigskip

\noindent\textbf{Theorem 4.3.} \textit{Let }$(\mathfrak{A},\omega,\tau
,K)$\textit{ be a }$C^{\ast}$\textit{-dynamical system with }$\omega$\textit{
a trace and }$\left(  K,\Sigma,\mu\right)  $\textit{ a F\o lner semigroup. Let
}$A\in\mathfrak{A}$\textit{ with }$\omega(A)>0$. \textit{Take any }%
$m_{0},\ldots,m_{k}\in N\cup\{0\}$. \textit{Assume that }$g\mapsto
\omega\left(  \prod_{j=0}^{k}\tau_{g^{m_{j}}}(A)\right)  $\textit{and}
$g\mapsto\Vert\tau_{g^{m_{j}}}(A)-A\Vert_{\omega}$\textit{ are }$\Sigma
$\textit{-measurable on }$K$\textit{ for }$j=0,1,\ldots,k$\textit{. Then there
exists a F\o lner net }$(\Lambda_{\alpha})$\textit{ in }$\left(  K,\Sigma
,\mu\right)  $\textit{ such that}
\[
\liminf_{\alpha}\frac{1}{\mu(\Lambda_{\alpha})}\int_{\Lambda_{\alpha}}\left|
\omega\left(  \prod_{j=0}^{k}\tau_{g^{m_{j}}}(A)\right)  \right|  d\mu(g)>0.
\]

\noindent\textbf{Proof.} This follows from Proposition 3.6 and Corollary 4.2
since $E=\{g\in K:\Vert\tau_{g^{m_{j}}}(A)-A\Vert_{\omega}<\varepsilon\text{
for }j=0,...,k\}$ is $\Sigma$-measurable. $\square$

\bigskip

Note that if for example $K$ is a topological semigroup and we assume that
$g\mapsto\tau_{g}(A)$ is continuous in $\mathfrak{A}$'s norm, then both
$g\mapsto\omega\left(  \prod_{j=0}^{k}\tau_{g^{m_{j}}}(A)\right)  $ and
$g\mapsto\Vert\tau_{g^{m_{i}}}(A)-A\Vert_{\omega}$ are continuous and hence
Borel measurable.

The Szemer\'{e}di property for a compact measure preserving dynamical system
$(X,\Sigma,\nu,T)$ with evolution over $\mathbb{N}$ is a special case of this
theorem, but note that $T$ need not be invertible in this case. Just let
$\omega(f):=\int_{X}fd\nu$ and $\tau(f):=f\circ T$ for all $f\in
\mathfrak{A}:=B_{\infty}(\Sigma)$ (see Section 2), let $\tau_{n}=\tau^{n}$ for
$n\in\mathbb{N}$, set $\Lambda_{N}:=\left\{  1,...,N\right\}  $ for all
$N\in\mathbb{N}$, and let $A=f$ be a positive function in $B_{\infty}(\Sigma)$
which is not $\nu$-a.e. zero. Keep in mind that the conclusion of Proposition
3.6 (and hence that of Theorem 4.3) holds for this choice of $\left(
\Lambda_{N}\right)  $, as is well known. The condition $\left|  \left|
\tau(f)\right|  \right|  \leq\left|  \left|  f\right|  \right|  $ follows
directly from $\tau$'s definition, while $\left|  \left|  \tau(f)\right|
\right|  _{\omega}=\left|  \left|  f\right|  \right|  _{\omega}$ expresses the
fact that $T$ is measure preserving, namely $\nu\circ T^{-1}=\nu$ as set
functions on $\Sigma$. More specifically (1.1) is obtained for these
assumptions by taking $f$ to be the characteristic function $\chi_{V}$ of a
set $V\in\Sigma$ with $\nu(V)>0$ and setting $m_{j}=j$

Lastly we note that Corollary 4.2 already contains much of the interesting
information regarding recurrence in a compact system, but under more general
conditions than Theorem 4.3. In the measure theoretic case in the previous
paragraph the inequality (4.1) becomes%
\[
\nu\left(  T^{-m_{0}n}V\cap...\cap T^{-m_{k}n}V\right)  >a
\]
for all $n$ in a relatively dense set in $\mathbb{N}$, where we have again
taken $A=\chi_{V}$ with $\nu(V)>0$.

\bigskip

\noindent\textbf{Acknowledgment}

The third author thanks the National Research Foundation for financial support.

\end{document}